\documentclass[a4paper,12pt,reqno]{amsart}
\usepackage[T1]{fontenc}
\usepackage[utf8]{inputenc}
\usepackage{amsmath,amssymb,mathtools}
\usepackage{needspace}
\newcommand{\lead}[1]{\par\smallskip\noindent\leavevmode#1\par\smallskip}
\usepackage[hidelinks]{hyperref}
\hypersetup{
 pdftitle={There are no Riesz bases of exponentials in balls and triangles},
 pdfauthor={Joaquim Ortega-Cerdà}
}
\usepackage[defersetup=false]{microtype}

\newtheorem{theorem}{Theorem}[section]
\newtheorem{lemma}[theorem]{Lemma}
\newtheorem{proposition}[theorem]{Proposition}
\newtheorem{corollary}[theorem]{Corollary}
\theoremstyle{remark}
\newtheorem{remark}[theorem]{Remark}
\numberwithin{equation}{section}

\newcommand{\R}{\mathbb R}
\newcommand{\C}{\mathbb C}
\newcommand{\Hh}{\mathcal H}
\newcommand{\Ff}{\mathcal F}
\newcommand{\one}{\mathbf 1}
\DeclareMathOperator{\Ran}{Ran}
\DeclareMathOperator{\supp}{supp}
\DeclareMathOperator{\conv}{conv}
\DeclareMathOperator{\spanop}{span}
\DeclareMathOperator{\dist}{dist}
\newcommand{\cP}{\mathcal P}

\title[No Riesz bases of exponentials]{There are no Riesz bases of exponentials in balls and triangles}
\author{Joaquim Ortega-Cerdà}
\address{Universitat de Barcelona\\Centre de Recerca Matemàtica (CRM)}
\email{jortega@ub.edu}
\thanks{The author is supported by grant PID2024-160033NB-I00, funded by
MICIU/AEI/10.13039/501100011033 and FEDER, UE, and by the Generalitat de
Catalunya through grant 2024 ICREA 00142.}
\date{September 25, 2026 (Version 2)}

\begin{document}
\raggedbottom

\begin{abstract}
We prove that balls in $\R^n$, $n\geq2$, and nondegenerate planar
triangles admit no Riesz basis of exponentials. The ball result extends
to every bounded convex domain with $C^2$ boundary in dimensions at least
two. The proof uses stationary limits and an obstruction for orthogonal
projections. We also construct complete and minimal exponential systems
with real frequencies on every bounded convex domain, with all
biorthogonal functions supported in its John ellipsoid.
\end{abstract}

\maketitle

\section{Introduction}

Fix an integer $n\geq2$.  Let $\Omega\subset\R^n$ be an open ball and, for $\lambda\in\R^n$, write
\[
 e_\lambda(x)=e^{2\pi i\lambda\cdot x}.
\]
We ask whether there is a discrete set $\Lambda\subset\R^n$ for which
$\{e_\lambda\}_{\lambda\in\Lambda}$ is a Riesz basis of $L^2(\Omega)$.  Thus the
synthesis map
\[
 (c_\lambda)\longmapsto \sum_{\lambda\in\Lambda}c_\lambda e_\lambda
\]
is an isomorphism from $\ell^2(\Lambda)$ onto $L^2(\Omega)$; equivalently, the exponential system is complete and there
are constants $0<A\leq B<\infty$ such that
\begin{equation}\label{eq:riesz}
 A\sum_{\lambda\in\Lambda}|c_\lambda|^2
 \leq
 \left\|\sum_{\lambda\in\Lambda}c_\lambda e^{2\pi i\lambda\cdot x}\right\|_{L^2(\Omega)}^2
 \leq
 B\sum_{\lambda\in\Lambda}|c_\lambda|^2
\end{equation}
for every finitely supported family $(c_\lambda)$.

The existence of an exponential Riesz basis for Euclidean balls is a longstanding
folklore problem in harmonic analysis. The author learned about this problem from Yurii Lyubarskii in the late 90's.  Bruna, Massaneda and Ortega-Cerdà
\cite[Section~5, p.~350]{BrunaMassanedaOrtegaCerda} explicitly formulated
the conjecture that balls in dimensions greater than one have no such basis.
The planar case was still recorded as open in 2025
\cite[Section~1.2]{LevTselishchev}.

The main result is the following.

\begin{theorem}\label{thm:ball}
For every $n\geq2$, no ball of positive radius in $\R^n$ admits a Riesz basis of exponentials.
\end{theorem}

The same conclusion holds for every nonempty bounded convex open set
with $C^2$ boundary in $\R^n$, $n\geq2$; see Corollary~\ref{cor:convex} in
Section~\ref{sec:convex}.
The argument also excludes exponential Riesz bases on every nondegenerate
triangle in $\R^2$ (Theorem~\ref{thm:triangle}).  More generally, Remark~\ref{rem:general-boundary} states a
boundary-measure criterion under which the same projection obstruction
applies to other open sets.

The failure of exponential Riesz bases for balls, triangles and smooth
convex domains makes it natural to ask whether these sets still admit
complete and minimal exponential systems. Here minimality means that
removing any one exponential destroys completeness, without a uniform
stability requirement. Olevskii explicitly posed the general existence
question in his 2021 Fields Institute lecture
\cite[Section~6, p.~7]{OlevskiiFields}, for measurable sets
$S\subset\R$ of finite positive measure. In
Section~\ref{sec:complete-minimal} we prove the convex-domain analogue in
every dimension: every nonempty bounded convex open set admits a complete
and minimal system of exponentials with real frequencies
(Theorem~\ref{cm:thm:main}). No boundary regularity or symmetry is needed,
and all the biorthogonal functions can be supported in the John ellipsoid
of $\overline\Omega$. This result addresses the higher-dimensional convex
case; the question for arbitrary measurable subsets of the line is not
settled by this construction.

The orthogonal version of this problem was initiated by
Fuglede \cite{Fuglede}.  Lev and Matolcsi \cite{LevMatolcsi} settled
Fuglede's conjecture for convex domains in every dimension: a convex body
admits an orthogonal basis of exponentials if and only if it tiles $\R^n$
by translations.  Iosevich, Katz and Pedersen
\cite{IosevichKatzPedersen} proved that balls in dimensions at least two
admit no orthogonal basis of exponentials.  Quantitative restrictions on
orthogonal exponential systems on a planar ball were obtained by Iosevich and
Kolountzakis \cite{IosevichKolountzakis}, and subsequently strengthened by
Zakharov \cite{Zakharov}.  These results concern pairwise orthogonality.
For Riesz bases, Kozma, Nitzan and Olevskii \cite{KozmaNitzanOlevskii}
constructed a bounded subset of the real line without an exponential
Riesz basis and obtained a quantitative lower bound for the Riesz constant
of a hypothetical basis on a planar ball, without excluding its existence.

Positive results are known for polygonal domains.
Lyubarskii and Rashkovskii \cite{LyubarskiiRashkovskii} constructed
exponential Riesz bases for centrally symmetric convex polygons whose
vertices lie on a lattice
using entire functions of two complex variables.
Debernardi and Lev \cite{DebernardiLev} later removed the arithmetic restriction, proving existence
for every centrally symmetric convex polygon and, more generally, for
centrally symmetric convex polytopes whose faces of all dimensions are
centrally symmetric.

A related construction principle is lattice multi-tiling: translates of a
bounded measurable domain by a lattice cover almost every point of
$\R^n$ the same finite number of times.  Grepstad and Lev \cite{GrepstadLev} proved
that such domains admit exponential Riesz bases under a Riemann
measurability assumption; Kolountzakis \cite{Kolountzakis} gave a simpler
proof and removed this regularity assumption.  Kozma and Nitzan proved
existence for finite unions of intervals \cite{KozmaNitzanIntervals}
and for finite unions of rectangles with sides parallel to the coordinate
axes \cite{KozmaNitzanRectangles}.

There is a function-theoretic reformulation of the problem.  If
\[
 PW_\Omega
 =
 \left\{
 F(z)=\int_\Omega f(x)e^{-2\pi iz\cdot x}\,dx : f\in L^2(\Omega)
 \right\},
 \qquad z\in\C^n,
\]
with the norm transported from $L^2(\Omega)$, then $PW_\Omega$ is a Hilbert space of
entire functions.  For real $\lambda$,
\[
 F(\lambda)=\langle f,e_\lambda\rangle_{L^2(\Omega)}.
\]
Hence $\Lambda$ is the frequency set of an exponential Riesz basis precisely
when the restriction map
\[
 F\longmapsto \bigl(F(\lambda)\bigr)_{\lambda\in\Lambda}
\]
is an isomorphism from $PW_\Omega$ onto $\ell^2(\Lambda)$.  Under the Fourier--Laplace identification, this restriction map is the
adjoint of the synthesis operator. Thus the question is whether this
Paley--Wiener space in $n$ variables has a complete interpolating sequence.

The proof relies on a geometric
feature of the ball: a nontrivial translate of its boundary meets the original
boundary in a set of $(n-1)$-dimensional surface measure zero.

The proof has four steps. First, separation of the frequency set allows
us to replace its points by orthonormal translates of a small smooth bump.
The Fourier transform of the bump is bounded away from zero on the domain,
so the Riesz basis property yields an estimate
$\|P-Q_\Lambda\|<1$ between two orthogonal projections on $L^2(\R^n)$.
Second, we take weak limits of translates of $\Lambda$ and construct an
invariant probability measure on the resulting compact space. This use of
weak limits originates in Beurling's work on sampling and interpolation
\cite{Beurling}; a modern account of the topology is given in
\cite[Section~4]{GOR}. Third, scalar Bochner measures and continuous box
averages yield spectral cutoffs $\Pi_t$ and comparison projections $M_t$
on the associated stationary Hilbert space, with
$\|\Pi_t-M_t\|\leq\gamma<1$ for almost every $t$. We obtain this estimate
by applying the original projection bound to compactly supported functions
and averaging the resulting scalar inner products. Finally, crossing the
boundary produces two strong limits of the cutoff projections with
strictly nested ranges. Norm continuity of $M_t$ places both limits at
distance at most $\gamma$ from $M_{t_0}$, giving a contradiction.

The translation-hull construction belongs to the established theory of
point set dynamical systems; see Lenz--Stollmann \cite{LenzStollmann} and
Bj\"orklund, Hartnick and Pogorzelski \cite{BjorklundHartnickPogorzelski}.
Its applications to sampling and interpolation are developed by
Enstad--Raum \cite{EnstadRaum}. Averaging inequalities over invariant
probability measures on hulls also plays a central role in
\cite{EnstadRaum,EnstadVanVelthoven}. Here we use this framework to
transfer a uniform distance estimate between projections to a stationary
Hilbert space, where a boundary crossing produces the obstruction.

Sections~\ref{sec:bumps}--\ref{sec:distance} establish the analytic
construction for any bounded open domain with boundary of Lebesgue measure
zero. Section~\ref{sec:sphere} supplies the geometric step for balls,
Section~\ref{sec:triangles} treats triangles, and Section~\ref{sec:convex}
extends the result to nonempty bounded convex open sets with $C^2$
boundary. Section~\ref{sec:complete-minimal} proves the existence of
complete minimal exponential systems on all bounded convex domains by
a separate Fourier-analytic construction.

\subsection*{Note} The non-existence of Riesz bases of exponentials in
balls and triangles has been independently obtained with a different
proof by Zhiqiang Wan in \cite{Wan}. We became aware of his result in the
final stage of preparing our manuscript. Since the methods are different
and our result also covers the full class of nonempty bounded convex open
sets with $C^2$ boundary, a class not treated in \cite{Wan}, we offer this
proof here.

\subsection*{Acknowledgments}
The author thanks Joaquim Bruna for extensive discussions on the manuscript.

\subsection*{Use of artificial intelligence}
A large language model was used for mathematical discussions and
assistance with the Lean formalization. The work began with an unsuccessful
attempt by the author to prove nonexistence of an exponential Riesz basis
on the disk, following a strategy inspired by
\cite{MateiMeyerOrtegaCerda}. During these discussions, the model corrected
the failed approach, leading to a proof based on representation theory.
The author subsequently simplified that proof into the present form.

The Lean formalization of the nonexistence results in
Sections~\ref{sec:bumps}--\ref{sec:convex} and the existence theorem in
Section~\ref{sec:complete-minimal} is available at the GitHub repository:
\href{https://github.com/joaquimortega/RieszEuclidean}{RieszEuclidean}.
The repository contains the Lean sources, a standalone
proof file, and verification records identifying the statements checked by
Comparator and Lean's kernel. The existence theorem retains the
Titchmarsh--Lions convolution support theorem as an explicit hypothesis.
\section{From complete interpolation to bump projections}\label{sec:bumps}

\lead{This section converts the Riesz basis assumption into a geometric
statement about two closed subspaces of $L^2(\R^n)$. Disjoint smooth bumps
model the coefficient space, and their projections onto the Fourier
subspace reproduce the original synthesis operator up to an invertible
multiplier.}

Throughout Sections~\ref{sec:bumps}--\ref{sec:distance}, $\Omega$ is a
bounded open subset of $\R^n$ of positive measure. Inner products are
linear in the first variable. We use the Fourier transform with the convention
\[
 (\Ff a)(x)=\int_{\R^n}a(\xi)e^{2\pi i\xi\cdot x}\,d\xi,
\]
initially defined on $L^1\cap L^2$, and set
\begin{equation}\label{eq:P}
 P=\Ff^{-1}M_{\one_\Omega}\Ff,\qquad
 p(y)=\int_\Omega e^{-2\pi iy\cdot x}\,dx.
\end{equation}
Thus $P$ is an orthogonal projection on $L^2(\R^n)$ that commutes with
translations. For compactly supported $a\in L^2(\R^n)$, Fourier inversion
and Fubini's theorem give
$(Pa)(v)=\int_{\R^n}p(v-w)a(w)\,dw$ almost everywhere.

\subsection{Separated frequencies and an exact bump representation}

\lead{The lower Riesz bound first separates the frequencies. We then choose
one sufficiently small bump so that its translates by the given frequencies are orthonormal and its
Fourier transform has no zeros on $\overline\Omega$. This gives an equivalence between the Riesz basis property and a
distance bound for projections, using the elementary criterion proved
in Lemma~\ref{lem:projection-distance} below.}

\begin{lemma}\label{lem:bumps}
The system $\{e_\lambda:\lambda\in\Lambda\}$ is a Riesz basis of
$L^2(\Omega)$ if and only if there are $\delta,r>0$, a real nonnegative
$b\in C_c^\infty(\R^n)$ and $0\leq\gamma<1$ such that $\Lambda$ is
$\delta$-separated,
\begin{equation}\label{eq:bump}
 \|b\|_2=1,\quad \supp b\subset B(0,r),\quad 2r<\delta,
 \quad \inf_{x\in\overline\Omega}|\Ff b(x)|>0,
\end{equation}
and the orthogonal projection $Q_\Lambda$ onto the closed span of
$\{b(\cdot-\lambda):\lambda\in\Lambda\}$ satisfies
\begin{equation}\label{eq:initial-distance-bound}
 \|P-Q_\Lambda\|\leq\gamma<1.
\end{equation}
\end{lemma}

\begin{proof}
Suppose first that the system is a Riesz basis. Apply the lower Riesz bound to coefficients $1$ and $-1$ at distinct
frequencies. Since $|e^{iu}-e^{iv}|\leq|u-v|$,
\[
 2A\leq\|e_\lambda-e_\eta\|_{L^2(\Omega)}^2
 \leq4\pi^2|\lambda-\eta|^2\int_\Omega|x|^2\,dx.
\]
The integral is finite and positive. This gives a positive separation
constant $\delta$. Comparing the volumes of disjoint balls of radius
$\delta/2$ also gives a uniform bound on the number of frequencies in any
ball of fixed radius.

Choose a real nonnegative smooth compactly supported $b_0$ with
$\|b_0\|_2=1$. For small $\varepsilon>0$ put
$b(\xi)=\varepsilon^{-n/2}b_0(\xi/\varepsilon)$. Its support can be made arbitrarily small, and
\[
 (\Ff b)(x)=\varepsilon^{n/2}(\Ff b_0)(\varepsilon x).
\]
Since $(\Ff b_0)(0)>0$ and $\overline\Omega$ is compact, this proves
\eqref{eq:bump} for a sufficiently small fixed $\varepsilon$.

For either direction, the conditions on $b$ and $\Lambda$ imply that the
translates $b_\lambda=b(\cdot-\lambda)$ have disjoint supports and
norm one. Consequently
\[
 W_\Lambda:\ell^2(\Lambda)\longrightarrow\Ran Q_\Lambda,
 \qquad W_\Lambda c=\sum_{\lambda\in\Lambda}c_\lambda b_\lambda,
\]
is unitary. For finitely supported $c$,
\begin{equation}\label{eq:bump-synthesis}
 \Ff(PW_\Lambda c)
 =\one_\Omega(\Ff b)\sum_{\lambda\in\Lambda}c_\lambda e_\lambda.
\end{equation}
Multiplication by $\Ff b$ is an isomorphism of $L^2(\Omega)$.
If the exponential system is a Riesz basis, its synthesis operator is
also an isomorphism onto $L^2(\Omega)$. Thus \eqref{eq:bump-synthesis}
extends to $\ell^2(\Lambda)$ and shows that
$P:\Ran Q_\Lambda\to\Ran P$ is an isomorphism.
Lemma~\ref{lem:projection-distance} below gives
\eqref{eq:initial-distance-bound}.

Conversely, suppose the stated separation, bump and distance conditions
hold. The same lemma shows that
$P:\Ran Q_\Lambda\to\Ran P$ is an isomorphism. Hence
\[
 c\longmapsto (\Ff b)^{-1}\Ff(PW_\Lambda c)|_\Omega
\]
is a bounded isomorphism from $\ell^2(\Lambda)$ onto $L^2(\Omega)$.
By \eqref{eq:bump-synthesis}, it agrees with exponential synthesis on
finitely supported sequences. The exponential system is therefore a
Riesz basis.
\end{proof}

\subsection{Distance between orthogonal projections}

\lead{We need two elementary facts about orthogonal projections. A distance
strictly smaller than one is equivalent to invertibility of the projection
between the ranges; consequently, two strictly nested ranges cannot both
be this close to the range of a third projection.}

For orthogonal projections $P$ and $Q$ on the same Hilbert space, we use
the operator-norm distance $\|P-Q\|$.

\Needspace{8\baselineskip}
\begin{lemma}\label{lem:projection-distance}
Let $P$ and $Q$ be orthogonal projections on a Hilbert space.  The map
\[
 P:\Ran Q\longrightarrow\Ran P
\]
is a bounded isomorphism if and only if
\begin{equation}\label{eq:projection-distance}
 \|P-Q\|<1.
\end{equation}
\end{lemma}

\begin{proof}
Set $A=P|_{\Ran Q}:\Ran Q\to\Ran P$. Its adjoint is
$A^*=Q|_{\Ran P}$. We first record the identity
\[
 \|P-Q\|
 =\max\{\|(I-P)Q\|,\|(I-Q)P\|\}.
\]
Indeed, for every $x$,
\[
 (P-Q)x=P(I-Q)x-(I-P)Qx.
\]
The two terms on the right are orthogonal, as are $(I-Q)x$ and $Qx$.
Since $\|P(I-Q)\|=\|(I-Q)P\|$, Pythagoras gives the upper bound in
the identity. The reverse bound follows by multiplying $P-Q$ on the left
by $I-P$ and $I-Q$.

Suppose first that $\|P-Q\|=\gamma<1$. For $u\in\Ran Q$ and
$y\in\Ran P$, Pythagoras gives
\[
 \|Au\|^2\geq(1-\gamma^2)\|u\|^2,
 \qquad
 \|A^*y\|^2\geq(1-\gamma^2)\|y\|^2.
\]
Thus $A$ is injective with closed range, while
$(\Ran A)^\perp=\ker A^*=\{0\}$ makes its range dense in $\Ran P$.
Hence $A$ is onto and has a bounded inverse.

Conversely, suppose $A$ is an isomorphism. If both ranges are zero, the
conclusion is immediate. Otherwise, put
$c=\|A^{-1}\|^{-1}\in(0,1]$. Since
$(A^*)^{-1}=(A^{-1})^*$, both $A$ and $A^*$ are bounded below by $c$.
Pythagoras therefore yields
\[
 \|(I-P)Q\|\leq\sqrt{1-c^2},
 \qquad
 \|(I-Q)P\|\leq\sqrt{1-c^2}.
\]
The identity above gives $\|P-Q\|\leq\sqrt{1-c^2}<1$.
\end{proof}

The following consequence will give the final contradiction.

\begin{lemma}\label{lem:inclusion}
Let $R_-$ and $R_+$ be orthogonal projections with
\[
 \Ran R_-\subsetneq\Ran R_+.
\]
There is no orthogonal projection $M$ such that
\[
 \|R_--M\|<1,
 \qquad
 \|R_+-M\|<1.
\]
\end{lemma}

\begin{proof}
By Lemma~\ref{lem:projection-distance}, $M$ maps each of the two ranges isomorphically onto
$\Ran M$.  Choose $w\in\Ran R_+\setminus\Ran R_-$.  Surjectivity on the
smaller range gives $u\in\Ran R_-$ with $Mu=Mw$.  Then
$0\neq w-u\in\Ran R_+$ and $M(w-u)=0$, contradicting injectivity on
$\Ran R_+$.
\end{proof}

Fix $b,\delta,r$ and $\gamma$ from Lemma~\ref{lem:bumps}. For every $\delta$-separated
set $\Gamma$, let $Q_\Gamma$ be the orthogonal projection onto the closed
span of the corresponding bumps. Explicitly,
\begin{equation}\label{eq:Q}
 Q_\Gamma a=\sum_{\lambda\in\Gamma}\langle a,b_\lambda\rangle b_\lambda.
\end{equation}

\section{Beurling weak limits and a stationary system}\label{sec:weak}

\lead{We enlarge the original frequency set to the compact space of all
its translation limits. The projection estimate persists throughout this
space. Averaging translates then provides an invariant probability measure
and a Hilbert space with a nonzero translation-invariant vector.}

\subsection{Weak convergence of separated sets}

\lead{Uniform separation prevents points from colliding in a limit. It
also makes local convergence of configurations equivalent to vague
convergence of their counting measures, which provides a convenient route
to compactness.}

We use Beurling's notion of weak convergence of point sets
\cite{Beurling}. For a modern treatment see \cite[Section~4]{GOR}.
A sequence of closed sets $\Gamma_j$ converges weakly to $\Gamma$ if,
for every $R,\varepsilon>0$, eventually
\begin{equation}\label{eq:weak-sets}
 \Gamma_j\cap B(0,R)\subset\Gamma+B(0,\varepsilon),\qquad
 \Gamma\cap B(0,R)\subset\Gamma_j+B(0,\varepsilon).
\end{equation}
For sequences with a common positive separation constant, this is equivalent
to vague convergence of their counting measures
$\eta_\Gamma=\sum_{\lambda\in\Gamma}\delta_\lambda$.
To see this, work in bounded regions whose boundaries contain no points
of the limit set. Weak convergence matches the finitely many points in
each such region, so continuous compactly supported test functions give
vague convergence. Conversely, testing on small disjoint balls around
limit points forces the matching points to exist. A sequence of unmatched
points in a fixed compact set would have a convergent subsequence and
would be detected by another nonnegative continuous test function.
Uniform separation prevents more than one point from approaching a limit
point.

The topology in \eqref{eq:weak-sets} is the restriction of the
Chabauty--Fell topology; see
\cite[Appendix~A, Proposition~A.1]{BjorklundHartnickPogorzelski}.
Lemma~\ref{lem:hull} is a standard consequence of compactness and
metrizability of that space and continuity of its translation action:
\cite[Proposition~2.2(ii)]{EnstadRaum}, with $U=B(0,\delta/2)$, shows
that the class of $\delta$-separated sets is closed. See also
\cite[Theorem~1.2, Proposition~1.3 and Remark~1.5(c)]{LenzStollmann}
for the Euclidean topology and hull construction. We include a proof
using counting measures.

\begin{lemma}\label{lem:hull}
The space of $\delta$-separated subsets of $\R^n$, including the empty
set, is compact and metrizable in the topology \eqref{eq:weak-sets}.
The translation action is jointly continuous in the translation parameter
and the configuration.
\end{lemma}

\begin{proof}
The separation bound gives a uniform bound on the mass of the counting
measures on each compact set. A diagonal extraction against a countable
dense collection of continuous compactly supported functions gives a
vaguely convergent subsequence, with a locally finite positive limit
measure $\eta$.

We verify that $\eta$ is again a counting measure. If
$x\in\supp\eta$, choose arbitrarily small balls $B(x,\varepsilon)$ of
diameter less than $\delta$ whose boundaries have $\eta$-measure zero.
Vague convergence gives
\[
 \eta_{\Gamma_j}(B(x,\varepsilon))
 \longrightarrow\eta(B(x,\varepsilon)).
\]
The left-hand side takes values in $\{0,1\}$, while the limit is positive
because $x\in\supp\eta$. Hence every such ball has mass one. Letting
$\varepsilon\downarrow0$ gives $\eta(\{x\})=1$. Thus $\eta$ is the
counting measure of its support. Two distinct support points at distance
less than $\delta$ would, for large $j$, force two points of $\Gamma_j$
to be less than $\delta$ apart. The support is therefore
$\delta$-separated; the zero measure corresponds to the empty set.

The same countable test functions metrize vague convergence on this
uniformly locally bounded family. Sequential compactness therefore gives
compactness. Finally, if $\Gamma_j\to\Gamma$ and $y_j\to y$, the local
matching criterion gives $\Gamma_j-y_j\to\Gamma-y$.
\end{proof}

Define the translation hull
\[
 X=\overline{\{\Lambda-y:y\in\R^n\}},\qquad T_y\Gamma=\Gamma-y.
\]
This is a compact metrizable translation-invariant space. Its elements
are precisely the weak limits of sequences of translates of $\Lambda$.

\subsection{Passing the projection estimate to weak limits}

\lead{Only finitely many bumps act on a compactly supported test function,
so local convergence of the centers gives strong convergence of the bump
projections. We also record a basic fact about strong limits that will
avoid repeating the same argument in the later constructions.}

A sequence of operators $A_j$ converges \emph{strongly} to $A$ if
$\|A_jf-Af\|\to0$ for every fixed vector $f$. We use the same
terminology for operator families depending on a real parameter.
Strong convergence is preserved under products of uniformly bounded
families: if $A_j\to A$ and $B_j\to B$ strongly, then
\begin{equation}\label{eq:strong-products}
 \|(A_jB_j-AB)f\|
 \leq\|A_j\|\,\|(B_j-B)f\|+\|(A_j-A)Bf\|
 \longrightarrow0.
\end{equation}
Self-adjointness passes to strong limits through inner products.
Consequently, a strong limit of orthogonal projections is an orthogonal
projection: the projections are uniformly bounded, and
\eqref{eq:strong-products} preserves idempotence.

The local matching of centers used in the next proof is the separated
case of \cite[Proposition~2.3]{EnstadRaum}; see also
\cite[Lemma~2.2]{EnstadVanVelthoven}. Compact support of the bumps
then gives strong convergence of the associated projections.

\begin{lemma}\label{lem:strong-Q}
If $\Gamma_j\to\Gamma$ in $X$, then $Q_{\Gamma_j}\to Q_\Gamma$
strongly on $L^2(\R^n)$. In particular,
\begin{equation}\label{eq:hull-distance-bound}
 \|P-Q_\Gamma\|\leq\gamma\qquad(\Gamma\in X).
\end{equation}
\end{lemma}

\begin{proof}
First take a compactly supported $a\in L^2(\R^n)$. Choose a ball
containing the closed $r$-neighborhood of its support in its interior
and having no point of $\Gamma$ on its boundary. Within this ball,
weak convergence and separation give, for large $j$, a bijection
between the finitely many centers of $\Gamma_j$ and $\Gamma$, with
matched centers converging. Only centers in this ball can contribute
to \eqref{eq:Q}. Translations of $b$ are continuous in $L^2$, so each
term of this finite sum converges in $L^2$. Thus
$Q_{\Gamma_j}a\to Q_\Gamma a$. Compactly supported functions are dense,
and all these projections have norm at most one, proving strong
convergence.

Conjugating by translations preserves $P$ and sends $Q_\Lambda$ to
the projection for the translated set. Thus every translate satisfies
\eqref{eq:initial-distance-bound}. Taking a strong limit and testing
on an arbitrary vector gives \eqref{eq:hull-distance-bound}.
\end{proof}

Since $P$ is nonzero, the empty set cannot belong to $X$: its projection
is zero and would have distance one from $P$. By
Lemma~\ref{lem:projection-distance} and \eqref{eq:bump-synthesis}, each
configuration in $X$ is again the frequency set of an exponential
Riesz basis for $L^2(\Omega)$. Preservation of the Riesz and frame
properties under hull limits is also established in
\cite[Proposition~3.8 and Theorem~3.9]{EnstadRaum}; for separated
sets the bounds are preserved without loss. Only the uniform projection
estimate will be needed below.

\subsection{Translation averaging}

\lead{Averaging the orbit over expanding boxes produces an invariant
probability measure. The associated translation operators are strongly
continuous, and the constant function has norm one. This invariant vector
will detect the spectral jump at the boundary.}

The following construction is the usual F\o lner averaging argument for
an invariant probability measure on a compact hull. In the sampling
setting, see \cite[proof of Theorem~2.16, Steps~2--3]{EnstadRaum};
see also the amenability discussion in
\cite[Section~2.3]{EnstadVanVelthoven}. Neither unique ergodicity nor
finite local complexity is needed here.

Put $C_R=[-R,R]^n$ and define a probability measure on $X$ by
\[
 \mu_R=\frac{1}{|C_R|}\int_{C_R}\delta_{T_v\Lambda}\,dv.
\]
Compactness gives a weakly convergent subsequence as $R\to\infty$;
write $\mu$ for one such limit. No uniqueness is needed. For
$F\in C(X)$ and fixed $y$,
\[
 \left|\int F\circ T_y\,d\mu_R-\int F\,d\mu_R\right|
 \leq\|F\|_\infty\frac{|(C_R+y)\mathbin\triangle C_R|}{|C_R|}
 \longrightarrow0.
\]
Passing to the limit shows that $\mu$ is invariant under every $T_y$.
Let
\begin{equation}\label{eq:U}
 \Hh=L^2(X,\mu),\qquad (U_yf)(\Gamma)=f(\Gamma+y).
\end{equation}
Thus $U_yf=f\circ T_{-y}$, and the operators $U_y$ form a unitary
representation of $\R^n$. Joint continuity of the action on the compact
space $X$ makes $y\mapsto U_yf$ continuous in the uniform norm for
$f\in C(X)$. Density of $C(X)$ in $\Hh$, together with unitarity,
extends strong continuity to every $f\in\Hh$. The Hilbert space $\Hh$
is separable because $X$ is compact and metrizable. The constant function
$\one$ has norm one and is fixed by every $U_y$.

\section{Bochner measures and stationary Fourier projections}\label{sec:bochner}

\lead{We now construct the stationary analogue of the Fourier projection
$P$. Scalar Bochner measures describe the frequencies of each vector, and
Fej\'er averages approximate the indicator of $\Omega$. A single control
measure makes the exceptional set of translation parameters independent
of the vector.}

From now on, assume in addition that $|\partial\Omega|=0$. We seek
projections on $\Hh$ that retain the stationary frequencies $\theta$
for which $t+\theta\in\Omega$. Both the parameter $t$ and the spectral
variable $\theta$ belong to $\R^n$.

\subsection{Scalar Fourier measures}

\lead{Bochner's theorem associates a positive measure to the correlation
function of each vector. The resulting norm identity turns convolution
kernels into spectral multipliers; it is all the Fourier calculus needed
below.}

For $f\in\Hh$, the continuous function
$y\mapsto\langle U_yf,f\rangle$ is positive definite. Bochner's theorem
\cite[Section~1.4.3, p.~19]{Rudin} gives a unique finite positive
Borel measure $\sigma_f$ on $\R^n$ such that
\begin{equation}\label{eq:bochner}
 \langle U_yf,f\rangle
 =\int_{\R^n}e^{2\pi iy\cdot\theta}\,d\sigma_f(\theta),
 \qquad \sigma_f(\R^n)=\|f\|^2.
\end{equation}
Uniqueness and the parallelogram identity give
\begin{equation}\label{eq:spectral-parallelogram}
 \sigma_{f+g}+\sigma_{f-g}=2\sigma_f+2\sigma_g,
 \qquad \sigma_{cf}=|c|^2\sigma_f.
\end{equation}
Since $U_y\one=\one$, uniqueness also gives
\begin{equation}\label{eq:constant-spectrum}
 \sigma_{\one}=\delta_0.
\end{equation}
This atom at zero will ensure that the boundary crossing produces a
nonzero change in the projection.

For $a\in L^1(\R^n)$ define
\[
 T_af=\int_{\R^n}a(y)U_yf\,dy,\qquad
 \phi_a(\theta)=\int_{\R^n}a(y)e^{2\pi iy\cdot\theta}\,dy.
\]
We call $a$ the \emph{kernel} of $T_a$ in translation coordinates and
$\phi_a$ its \emph{symbol} in frequency coordinates. The integral defining
$T_af$ is a Bochner integral in $\Hh$, with norm at most
$\|a\|_1\|f\|$. Expanding its squared norm and using
\eqref{eq:bochner} gives
\begin{align*}
 \|T_af\|^2
 &=\int\!\!\int a(y)\overline{a(w)}
       \langle U_{y-w}f,f\rangle\,dy\,dw\\
 &=\int_{\R^n}|\phi_a(\theta)|^2\,d\sigma_f(\theta).
\end{align*}
Fubini applies because $a\in L^1$ and $\sigma_f$ is finite.
The group law and unitarity give, for $a,c\in L^1(\R^n)$,
\begin{equation}\label{eq:multiplier-norm}
 \|T_af\|^2=\int|\phi_a|^2\,d\sigma_f,\qquad
 T_aT_c=T_{a*c},\qquad T_a^*=T_{a^*},
\end{equation}
where $a^*(y)=\overline{a(-y)}$. In particular
$\|T_a\|\leq\|\phi_a\|_\infty$. Convolution of kernels corresponds
to multiplication of symbols, and $\phi_{a^*}=\overline{\phi_a}$.
If two kernels have the same symbol, their operators agree by the
norm identity. These facts provide the scalar calculus used below.

\subsection{One measure to control all exceptional sets}

\lead{Different vectors may have different spectral measures. A weighted
sum over a countable dense set gives one measure whose null sets are null
for every vector, allowing us to choose the later boundary point only
once.}

Choose a dense sequence $(h_j)_{j\geq1}$ in the unit sphere of $\Hh$
and set
\begin{equation}\label{eq:control}
 \sigma=\sum_{j=1}^\infty2^{-j}\sigma_{h_j}.
\end{equation}
This is a positive measure of mass one. For every Borel set $E$,
\begin{equation}\label{eq:control-null}
 \sigma(E)=0\quad\Longrightarrow\quad
 \sigma_f(E)=0\quad(f\in\Hh).
\end{equation}
Indeed, approximate $f$ by scalar multiples $g$ of the $h_j$.
Then $\sigma_g(E)=0$ and \eqref{eq:spectral-parallelogram}, applied to $f-g$ and $g$, gives
\[
 \sigma_f(E)\leq2\sigma_{f-g}(E)+2\sigma_g(E)
 \leq2\|f-g\|^2\longrightarrow0.
\]
Thus a single null-set condition applies simultaneously to all vectors
in $\Hh$.

\subsection{Continuous Fej\'er filters}

\lead{The overlap of a box with its translates generates a nonnegative
approximate identity. Its compactly supported inverse Fourier kernel
serves two purposes: it defines bounded stationary multipliers here and
reappears in the spatial averages of the next section.}

The normalized overlap of the boxes $C_R=[-R,R]^n$ is
\begin{equation}\label{eq:fejer-weights}
 w_R(y)=\frac{|C_R\cap(C_R-y)|}{|C_R|}
 =\prod_{j=1}^n\left(1-\frac{|y_j|}{2R}\right)_+.
\end{equation}
Put
\[
 K_R(x)=\frac{1}{|C_R|}\left|\int_{C_R}e^{2\pi iv\cdot x}\,dv\right|^2,
 \qquad m_R=K_R*\one_\Omega.
\]
Plancherel gives $K_R\geq0$ and $\int K_R=1$. The scaling identity
$K_R(x)=R^nK_1(Rx)$ shows that for each $\varepsilon>0$,
$\int_{|x|>\varepsilon}K_R(x)\,dx\to0$. Thus $(K_R)_{R>0}$ is an approximate
identity as $R\to\infty$. Since $\one_\Omega$ is locally constant away from its boundary,
\begin{equation}\label{eq:fejer-convergence}
 0\leq m_R\leq1,\qquad
 m_R(x)\longrightarrow\one_\Omega(x)\quad(x\notin\partial\Omega).
\end{equation}
Expanding the square in $K_R$ shows that $K_R=\Ff w_R$. Consequently,
with $p$ as in \eqref{eq:P},
\begin{equation}\label{eq:fejer-coefficients}
 m_R(x)=\int_{\R^n}w_R(y)p(y)e^{2\pi iy\cdot x}\,dy.
\end{equation}
This identity follows by Fubini: $w_R$ has compact support and $p$ is
bounded. The same overlap weight will appear in the averaged inner
products in Section~\ref{sec:distance}.

Define
\begin{equation}\label{eq:finite-Pi}
 \Pi_{t,R}f=\int_{\R^n}w_R(y)p(y)e^{2\pi it\cdot y}U_yf\,dy.
\end{equation}
This is $T_a$ for an integrable kernel with symbol $m_R(t+\theta)$.
By \eqref{eq:multiplier-norm} and \eqref{eq:fejer-convergence}, each
$\Pi_{t,R}$ is a self-adjoint contraction.

\subsection{Constructing the cutoff projections}

\lead{For almost every translation parameter, the translated boundary
has zero control measure. At these parameters, the approximate symbols
converge in every relevant $L^2$ space, producing orthogonal projections
and an exact formula for their differences.}

\begin{proposition}\label{prop:cutoffs}
The set
\begin{equation}\label{eq:good-t}
 \mathcal G=\{t\in\R^n:\sigma(\partial\Omega-t)=0\}
\end{equation}
is measurable and has full Lebesgue measure. For every $t\in\mathcal G$,
$\Pi_{t,R}$ converges strongly as $R\to\infty$ to an orthogonal
projection $\Pi_t$. For $s,t\in\mathcal G$ and $f\in\Hh$,
\begin{align}
 \|\Pi_tf\|^2
 &=\int_{\R^n}\one_\Omega(t+\theta)\,d\sigma_f(\theta),
 \label{eq:Pi-norm}\\
 \|(\Pi_t-\Pi_s)f\|^2
 &=\int_{\R^n}|\one_\Omega(t+\theta)-\one_\Omega(s+\theta)|^2
       \,d\sigma_f(\theta).
 \label{eq:Pi-difference}
\end{align}
\end{proposition}

\begin{proof}
The function $(t,\theta)\mapsto\one_{\partial\Omega}(t+\theta)$ is
Borel measurable, so its integral against $\sigma$ is measurable.
Tonelli's theorem and the assumption $|\partial\Omega|=0$ give
\begin{equation}\label{eq:good-ae}
 \int_{\R^n}\sigma(\partial\Omega-t)\,dt
 =\int_{\R^n}|\partial\Omega-\theta|\,d\sigma(\theta)=0.
\end{equation}
Thus $\mathcal G$ is conull. Fix $t\in\mathcal G$. By
\eqref{eq:control-null}, the boundary exception is null for every
$\sigma_f$. The difference norm formula for the integrable kernels gives
\[
 \|(\Pi_{t,R}-\Pi_{t,S})f\|^2
 =\int|m_R(t+\theta)-m_S(t+\theta)|^2\,d\sigma_f(\theta).
\]
The bounded symbols converge in $L^2(\sigma_f)$ by dominated convergence.
Hence $\Pi_{t,R}f$ is Cauchy and has a limit $\Pi_tf$. Uniform boundedness
makes $\Pi_t$ a bounded linear contraction. Similarly,
\[
 \|(\Pi_{t,R}^2-\Pi_{t,R})f\|^2
 =\int|m_R(t+\theta)^2-m_R(t+\theta)|^2\,d\sigma_f(\theta)
 \longrightarrow0.
\]
By \eqref{eq:strong-products}, $\Pi_{t,R}^2\to\Pi_t^2$ strongly,
so the preceding limit implies $\Pi_t^2=\Pi_t$. Self-adjointness passes
through inner products, so $\Pi_t$ is orthogonal. Passing to limits in the norm identities for
$\Pi_{t,R}$ and $\Pi_{t,R}-\Pi_{s,R}$ proves the asserted formulas.
\end{proof}

We call $\theta\mapsto\one_\Omega(t+\theta)$ the cutoff symbol of
$\Pi_t$. Here $t$ is a parameter and $\theta$ is a frequency of the
stationary translation action. Formula~\eqref{eq:Pi-norm} says that
$\|\Pi_tf\|^2$ is the spectral mass of $f$ in $\Omega-t$.
Formula~\eqref{eq:Pi-difference} will allow us to construct the two
boundary limits.

\section{Preserving the distance bound}\label{sec:distance}

\lead{We transfer the original distance estimate to the stationary Hilbert
space without losing its constant. First we construct a norm-continuous
family of comparison projections $M_t$. We then insert stationary vectors
into compactly supported functions on $\R^n$ and average their scalar
inner products.}

\subsection{The covariant bump kernel}

\lead{The bump projection has a continuous kernel that vanishes when its
two spatial arguments are sufficiently far apart. Covariance turns this
kernel into operators on $\Hh$; its projection identity and bounded
support give idempotence and norm continuity, respectively.}

The projection $Q_\Gamma$ has the continuous kernel
\begin{equation}\label{eq:q-kernel}
 q_\Gamma(v,w)=\sum_{\lambda\in\Gamma}
 b(v-\lambda)\overline{b(w-\lambda)}.
\end{equation}
For fixed $v$ at most one summand is nonzero. Thus, uniformly in $\Gamma$,
\begin{equation}\label{eq:q-bounds}
 |q_\Gamma(v,w)|\leq\|b\|_\infty^2,\qquad
 q_\Gamma(v,w)=0\quad\text{if }|v-w|>2r.
\end{equation}
Local finiteness, weak convergence of the configurations, and continuity
of $b$ imply joint continuity in $(\Gamma,v,w)$. The kernel satisfies
\begin{align}
 q_{\Gamma+y}(v,w)&=q_\Gamma(v-y,w-y),\label{eq:q-covariance}\\
 q_\Gamma(w,v)&=\overline{q_\Gamma(v,w)},\label{eq:q-hermitian}\\
 \int q_\Gamma(v,u)q_\Gamma(u,w)\,du&=q_\Gamma(v,w).
 \label{eq:q-projection}
\end{align}
The last identity follows by expanding the locally finite sums and using
orthonormality of the bumps. The integrals are absolutely convergent because $b$ has compact support.

Set $k_y(\Gamma)=q_\Gamma(0,-y)$, and write $B_y$ for multiplication by
$k_y$ on $\Hh$. The map $y\mapsto k_y$ is continuous with values in
$C(X)$: joint continuity is uniform on compact subsets of $X\times\R^n$.
In particular $y\mapsto B_yU_yf$ is continuous for every $f\in\Hh$.
Moreover it vanishes for $|y|>2r$. Define
\begin{equation}\label{eq:M}
 M_tf=\int_{\R^n}e^{2\pi it\cdot y}B_yU_yf\,dy.
\end{equation}
This is an $\Hh$-valued integral with norm at most
$\|f\|\int\|k_y\|_\infty\,dy$.

Covariant operator families over point set dynamical systems are studied
in \cite[Sections~2--3]{LenzStollmann}. The construction here uses the
explicit bump kernel to obtain operators on $L^2(X,\mu)$; the next
proof verifies their projection property and norm continuity directly.

\begin{proposition}\label{prop:M}
Every $M_t$ is an orthogonal projection. The family is Lipschitz continuous
in operator norm:
\begin{equation}\label{eq:M-continuity}
 \|M_t-M_s\|\leq L|t-s|,\qquad
 L=2\pi\int_{\R^n}|y|\,\|k_y\|_\infty\,dy<\infty.
\end{equation}
\end{proposition}

\begin{proof}
For bounded $f,g$, invariance of $\mu$ and the change of variables
$\Gamma'=\Gamma+y$ give
\begin{align*}
 \langle B_yU_yf,g\rangle
 &=\int_X q_\Gamma(0,-y)f(\Gamma+y)\overline{g(\Gamma)}\,d\mu(\Gamma)\\
 &=\int_X q_{\Gamma'}(y,0)f(\Gamma')
                  \overline{g(\Gamma'-y)}\,d\mu(\Gamma')\\
 &=\langle f,B_{-y}U_{-y}g\rangle.
\end{align*}
Density gives $(B_yU_y)^*=B_{-y}U_{-y}$. Taking the adjoint in
\eqref{eq:M} and changing $y$ to $-y$ therefore gives $M_t^*=M_t$.

For idempotence, covariance implies
$k_u(\Gamma+y)=q_\Gamma(-y,-y-u)$. First let $f$ be bounded.
Expanding $M_t^2$ and putting $v=y+u$ yields
\begin{align*}
 (M_t^2f)(\Gamma)
 &=\int e^{2\pi it\cdot v}
       \left[\int q_\Gamma(0,-y)q_\Gamma(-y,-v)\,dy\right]
       f(\Gamma+v)\,dv\\
 &=\int e^{2\pi it\cdot v}q_\Gamma(0,-v)f(\Gamma+v)\,dv
 =(M_tf)(\Gamma),
\end{align*}
where the second equality is \eqref{eq:q-projection}, with the change
of variable $u=-y$. The original double integral has $|y|,|u|\leq2r$
and uniformly bounded coefficients, so Fubini applies. The equality
extends to every $f\in\Hh$ by density and boundedness. Hence $M_t^2=M_t$,
and $M_t$ is an orthogonal projection.

Finally $|e^{2\pi it\cdot y}-e^{2\pi is\cdot y}|\leq2\pi|t-s||y|$.
Applying the triangle inequality to the vector integral proves
\eqref{eq:M-continuity}. Finiteness follows from \eqref{eq:q-bounds}.
\end{proof}

Norm continuity of $M_t$ comes from the integrable, compactly supported
kernel. The translation operators $U_y$ themselves need only be strongly
continuous. This distinction is essential at the boundary crossing.

\subsection{Compactly supported functions from a stationary vector}

\lead{We next restrict a stationary vector to a finite box along each
translation orbit and add a modulation. The normalization is chosen so
that averaging the squared spatial norm recovers exactly the original
Hilbert-space norm.}

Fix $R>0$, $t\in\R^n$, and abbreviate $C=C_R$. For $f\in\Hh$ define,
for almost every $\Gamma$,
\begin{equation}\label{eq:finite-test}
 a_f(\Gamma)(v)=|C|^{-1/2}\one_C(v)e^{-2\pi it\cdot v}f(\Gamma-v).
\end{equation}
These are compactly supported functions of $v\in\R^n$. Tonelli's theorem and
translation invariance give
\begin{equation}\label{eq:finite-energy}
 \int_X\|a_f(\Gamma)\|_2^2\,d\mu(\Gamma)
 =\frac{1}{|C|}\int_C\int_X|f(\Gamma-v)|^2\,d\mu(\Gamma)\,dv
 =\|f\|^2.
\end{equation}
In particular $a_f(\Gamma)\in L^2(\R^n)$ for almost every $\Gamma$.
Borel representatives of $f$ give joint measurability, and the same
identity applied to a representative that vanishes almost everywhere
shows that the resulting integrals do not depend on representatives.
Thus these averages are well defined despite the uncountable family of
translation parameters.

\subsection{Computing the two averaged inner products}

\lead{The two spatial projections can now be averaged in parallel. A
change to the difference variable produces the same box-overlap weight in
both calculations, identifying the averages with $M_{t,R}$ and
$\Pi_{t,R}$.}

Let
\begin{equation}\label{eq:finite-M}
 M_{t,R}f=\int w_R(y)e^{2\pi it\cdot y}B_yU_yf\,dy.
\end{equation}
For bounded $f,g$, expansion of the compact spatial integrals gives
\begin{align*}
 &\int_X\langle Q_\Gamma a_f(\Gamma),a_g(\Gamma)\rangle\,d\mu(\Gamma)\\
 &\quad=\frac{1}{|C|}\int_{C\times C}e^{2\pi it\cdot(v-w)}
   \int_X q_\Gamma(v,w)f(\Gamma-w)
                   \overline{g(\Gamma-v)}\,d\mu(\Gamma)\,dv\,dw.
\end{align*}
In the inner integral, make the change of variables $\Gamma'=\Gamma-v$. For $y=v-w$,
covariance turns the kernel into $q_{\Gamma'}(0,-y)$, and the two
functions become $f(\Gamma'+y)$ and $g(\Gamma')$. The inner integral
is therefore $\langle B_yU_yf,g\rangle$. For every integrable function
$H$ of the difference variable,
\[
 \frac{1}{|C|}\int_{C\times C}H(v-w)\,dv\,dw
 =\int w_R(y)H(y)\,dy.
\]
Consequently
\begin{equation}\label{eq:average-Q}
 \int_X\langle Q_\Gamma a_f,a_g\rangle\,d\mu
 =\langle M_{t,R}f,g\rangle.
\end{equation}

The kernel $p(v-w)$ of $P$ is independent of $\Gamma$. The same
calculation, now with inner integral $\langle U_{v-w}f,g\rangle$, gives
\begin{align}
 \int_X\langle Pa_f,a_g\rangle\,d\mu
 &=\int w_R(y)p(y)e^{2\pi it\cdot y}\langle U_yf,g\rangle\,dy\notag\\
 &=\langle\Pi_{t,R}f,g\rangle.
 \label{eq:average-P}
\end{align}
All kernels are bounded on $C\times C$. Thus the expansions for bounded
$f,g$ are justified by absolute integrability and Fubini's theorem. For general $f,g\in\Hh$,
Cauchy--Schwarz, contractivity of $P,Q_\Gamma$ and
\eqref{eq:finite-energy} bound each averaged form by $\|f\|\|g\|$.
The identities extend by density.

\subsection{Applying the distance bound and passing to the limit}

\lead{The averaged identities retain the original operator-norm bound
because the test functions preserve the averaged squared norm. We then
pass to a strong limit for the spectral cutoffs and a norm limit for the
comparison operators.}

Related arguments integrating frame inequalities over an invariant hull
measure occur in \cite[Theorem~3.12]{EnstadRaum} and
\cite[Theorem~3.5]{EnstadVanVelthoven}. Those results relate frame
bounds to covolumes. Here the averaged forms
\eqref{eq:average-Q}--\eqref{eq:average-P} instead transfer the
operator-norm distance bound between the two projections, with the same
constant $\gamma$.

\begin{proposition}\label{prop:distance}
For every $t\in\mathcal G$,
\begin{equation}\label{eq:stationary-distance-bound}
 \|\Pi_t-M_t\|\leq\gamma<1.
\end{equation}
\end{proposition}

\begin{proof}
Subtracting \eqref{eq:average-Q} from \eqref{eq:average-P} gives
\begin{equation}\label{eq:compression-identity}
 \langle(\Pi_{t,R}-M_{t,R})f,g\rangle
 =\int_X\langle(P-Q_\Gamma)a_f(\Gamma),a_g(\Gamma)\rangle\,d\mu.
\end{equation}
Apply \eqref{eq:hull-distance-bound} for almost every $\Gamma$, followed
by Cauchy--Schwarz in $X$ and \eqref{eq:finite-energy}:
\begin{align*}
 |\langle(\Pi_{t,R}-M_{t,R})f,g\rangle|
 &\leq\gamma\int_X\|a_f(\Gamma)\|_2\|a_g(\Gamma)\|_2\,d\mu\\
 &\leq\gamma\|f\|\|g\|.
\end{align*}
Thus $\|\Pi_{t,R}-M_{t,R}\|\leq\gamma$ for every $t$ and $R$.
By \eqref{eq:q-bounds} and dominated convergence,
\begin{equation}\label{eq:M-norm-limit}
 \|M_{t,R}-M_t\|
 \leq\int_{|y|\leq2r}|1-w_R(y)|\,\|k_y\|_\infty\,dy
 \longrightarrow0,
\end{equation}
uniformly in $t$. For $t\in\mathcal G$, Proposition~\ref{prop:cutoffs}
gives strong convergence $\Pi_{t,R}\to\Pi_t$. Testing the distance
bound on each vector and taking the limit proves
\eqref{eq:stationary-distance-bound}.
\end{proof}

The finite averages $M_{t,R}$ need not be projections. Their norm limit
$M_t$ is a projection by Proposition~\ref{prop:M}. At a boundary parameter,
the spectral cutoffs will have two different strong limits, whereas
$M_t$ has a single operator-norm limit. This gives a common comparison
projection for the two sides.

\section{Crossing the sphere}\label{sec:sphere}

\lead{The analytic construction is now complete. For a ball, we choose a
boundary point at which the only spectral boundary mass is at frequency
zero. Approaching this point from the two sides adds that component
without removing any other one, creating strictly nested limiting ranges.}

We now prove Theorem~\ref{thm:ball}. Suppose that a ball
$\Omega=B(a,\rho)\subset\R^n$, $n\geq2$, admits an exponential Riesz
basis. The preceding sections give projections $\Pi_t$ and $M_t$ with
\eqref{eq:stationary-distance-bound} for $t\in\mathcal G$. We will
construct two strong limits at a boundary point $t_0$, with strictly
nested ranges and both at distance at most $\gamma<1$ from $M_{t_0}$.

\subsection{Choosing a boundary point with no other spectral boundary mass}

\lead{A sphere has surface-measure-zero intersection with each of its
nontrivial translates. Tonelli's theorem reverses the roles of boundary
points and frequencies, yielding a point whose translated boundary has
no control-measure mass away from zero.}

Write $S=\partial\Omega$ and let $\nu=\mathcal H^{n-1}|_S$ be surface
measure on $S$. The measure $\nu$ is finite and nonzero.

\begin{lemma}\label{lem:sphere}
For every nonzero $\theta\in\R^n$,
\begin{equation}\label{eq:sphere-intersection}
 \nu\bigl(S\cap(S-\theta)\bigr)=0.
\end{equation}
Moreover, $S$ has Lebesgue measure zero.
\end{lemma}

\begin{proof}
If $x\in S\cap(S-\theta)$, then $|x-a|^2=\rho^2$ and
$|x+\theta-a|^2=\rho^2$. Subtraction gives
\[
 2(x-a)\cdot\theta=-|\theta|^2.
\]
For $\theta\ne0$, this is a hyperplane. Its intersection with the sphere
is empty, a point, or an $(n-2)$-dimensional sphere, hence has zero
$\mathcal H^{n-1}$-measure. Finally a sphere in $\R^n$ has zero
$n$-dimensional Lebesgue measure.
\end{proof}

We need the reverse viewpoint: fix a boundary point $t_0$ and control
the frequencies $\theta$ for which $t_0+\theta\in S$.
Tonelli's theorem interchanges these roles. Using the control measure $\sigma$
from \eqref{eq:control}, we obtain
\begin{equation}\label{eq:boundary-fubini}
 \begin{aligned}
 \int_{S}\sigma\bigl((S-t)\setminus\{0\}\bigr)\,d\nu(t)
 &=\int_{\R^n\setminus\{0\}}
     \nu\bigl(S\cap(S-\theta)\bigr)\,d\sigma(\theta)\\
 &=0.
 \end{aligned}
\end{equation}
The integrand on the left is nonnegative and $\nu(S)>0$.  Thus we can
choose $t_0\in S$ such that
\begin{equation}\label{eq:boundary-null}
 \sigma\bigl((S-t_0)\setminus\{0\}\bigr)=0.
\end{equation}
By \eqref{eq:control-null}, the same set is null for every $\sigma_f$.
This is why we introduced one control measure: the chosen point works
for every vector in $\Hh$ simultaneously.

Zero must be excluded from this null-set statement.  Indeed,
$\sigma_{\one}=\delta_0$, so the constant vector has spectral mass
there.  This is precisely the mass that will enter the cutoff.

\subsection{Approaching the point from outside and inside}

\lead{The cutoff projection need not be defined at the chosen boundary
point. We therefore approach it through good parameters in ambient open
sets on the two sides, rather than along a prescribed line. The limiting
symbols agree everywhere except at the zero-frequency atom, up to a common
spectral null set.}

The good set $\mathcal G$ has full Lebesgue measure in $\R^n$, but this
does not imply that it contains $t_0$.  In fact $t_0\notin\mathcal G$:
since $0\in S-t_0$ and $\sigma_{\one}=\delta_0$, condition
\eqref{eq:good-t} cannot hold there.  We therefore approach $t_0$
through good parameters and take limits of their projections.

Every neighborhood of a point on the sphere meets both $\Omega$ and
its strict exterior in nonempty open sets.  Each has positive Lebesgue
measure, so each meets $\mathcal G$.  We can consequently choose
\[
 t_j^-\in\mathcal G\cap(\R^n\setminus\overline\Omega),\qquad
 t_j^+\in\mathcal G\cap\Omega,\qquad t_j^\pm\longrightarrow t_0.
\]
These choices use open neighborhoods in the ambient space.  There is
no need to restrict the good parameters to a prescribed radial line.

Write the corresponding cutoff symbols as
\[
 a_j^\pm(\theta)=\one_\Omega(t_j^\pm+\theta).
\]
If $t_0+\theta$ is away from $S$, membership in $\Omega$ is locally
constant, so both symbols eventually equal
$\one_\Omega(t_0+\theta)$.  At $\theta=0$, the exterior symbol is
always zero and the interior symbol is always one.  All remaining
cases lie in $(S-t_0)\setminus\{0\}$, which is null for every
$\sigma_f$ by \eqref{eq:boundary-null} and \eqref{eq:control-null}.
Thus, for every $f\in\Hh$, we have $a_j^\pm\to m_\pm$
$\sigma_f$-almost everywhere, where
\begin{align}
 m_-(\theta)&=\one_{\Omega-t_0}(\theta),\label{eq:mminus}\\
 m_+(\theta)&=\one_{\Omega-t_0}(\theta)+\one_{\{0\}}(\theta).
 \label{eq:mplus}
\end{align}
Since $\Omega$ is open and $t_0\in S$, zero does not belong to
$\Omega-t_0$.  Both limits are therefore indicator functions, and
$m_-\leq m_+$.  The only added spectral component is at zero.

\subsection{Turning the limiting symbols into projections}

\lead{The exact difference formula converts convergence of the indicator
functions into strong convergence of the projections. Their limiting norm
identities will then allow us to read inclusion of the ranges directly
from inclusion of the spectral sets.}

For either choice of sign, dominated convergence gives
$a_j^\pm\to m_\pm$ in $L^2(\sigma_f)$ for every $f\in\Hh$.
Formula~\eqref{eq:Pi-difference} therefore yields
\[
 \|(\Pi_{t_j^\pm}-\Pi_{t_k^\pm})f\|
 =\|a_j^\pm-a_k^\pm\|_{L^2(\sigma_f)}
 \longrightarrow0.
\]
Thus $(\Pi_{t_j^\pm}f)_j$ is Cauchy for every $f$. Since the projections
are contractions, their limits define bounded linear operators $R_-,R_+$
with
\begin{equation}\label{eq:boundary-strong}
 \Pi_{t_j^-}\longrightarrow R_-,\qquad
 \Pi_{t_j^+}\longrightarrow R_+
 \quad\text{strongly}.
\end{equation}
By the strong-convergence facts recorded in Section~\ref{sec:weak},
$R_-$ and $R_+$ are orthogonal projections. Passing to the limit in
\eqref{eq:Pi-norm} gives
\begin{equation}\label{eq:boundary-norm}
 \|R_\pm f\|^2=\int_{\R^n}m_\pm(\theta)\,d\sigma_f(\theta).
\end{equation}
This formula lets us read the order of the limiting projections from
the order of their symbols.

\subsection{Strict inclusion of subspaces and the contradiction}

\lead{Order of the limiting symbols gives inclusion of the ranges, and
the constant vector makes this inclusion strict. Norm continuity of the
comparison family then places both limiting projections within distance
less than one of the same projection, contradicting the elementary
projection lemma.}

Since $m_-\leq m_+$, \eqref{eq:boundary-norm} implies
$\|R_-f\|\leq\|R_+f\|$ for every $f$.  In particular,
$R_+f=0$ implies $R_-f=0$.  Thus
\[
 \ker R_+\subset\ker R_-,\qquad
 \Ran R_-\subset\Ran R_+,
\]
where the second inclusion follows by taking orthogonal complements.

To see that this inclusion is strict, use the constant vector.  Its
spectral measure is $\delta_0$, so \eqref{eq:boundary-norm} gives
\[
 R_-\one=0,\qquad \|R_+\one\|=\|\one\|=1.
\]
By Pythagoras, an orthogonal projection preserves a vector's norm only
when it fixes that vector.  Hence $R_+\one=\one$: the constant
function belongs to the larger range and is orthogonal to the smaller
one.  We have proved
\begin{equation}\label{eq:strict-inclusion}
 \Ran R_-\subsetneq\Ran R_+.
\end{equation}
More precisely, the difference is an orthogonal projection, and
\begin{equation}\label{eq:zero-frequency-jump}
 \|(R_+-R_-)f\|^2=\sigma_f(\{0\}),\qquad
 (R_+-R_-)\one=\one.
\end{equation}
The jump need not have rank one; the constant vector only ensures that
it is nonzero.

Both limits also inherit a common distance bound. By
Proposition~\ref{prop:M}, $M_{t_j^\pm}\to M_{t_0}$ in operator norm.
For every $f\in\Hh$,
\begin{align*}
 \|(R_\pm-M_{t_0})f\|
 &=\lim_j\|(\Pi_{t_j^\pm}-M_{t_0})f\|\\
 &\leq\limsup_j\left(\gamma+\|M_{t_j^\pm}-M_{t_0}\|\right)\|f\|
 =\gamma\|f\|.
\end{align*}
Consequently
\begin{equation}\label{eq:boundary-distance}
 \|R_--M_{t_0}\|\leq\gamma<1,\qquad
 \|R_+-M_{t_0}\|\leq\gamma<1.
\end{equation}
This comparison is valid even though $t_0\notin\mathcal G$: $M_t$ is
defined and is a projection for every $t\in\R^n$.
Lemma~\ref{lem:inclusion} contradicts \eqref{eq:strict-inclusion} and
\eqref{eq:boundary-distance}. This proves Theorem~\ref{thm:ball}.

\begin{remark}[The role of the stationary space]\label{rem:stationary-space}
The invariant probability measure supplies a nonzero constant vector
with spectral measure $\delta_0$. The ordinary translation representation
$(V_ya)(v)=a(v-y)$ on $L^2(\R^n)$ has no nonzero invariant vector.
With our Fourier convention, its scalar spectral measures are
$|\Ff a(\theta)|^2\,d\theta$, so changing a cutoff at the single point
zero has no effect on that space. On $\Hh$, the same change gives
$R_-\one=0$ and $R_+\one=\one$. The averaging argument brings both this
invariant vector and the original distance bound into the same stationary
space.
\end{remark}

\begin{remark}\label{rem:interval}
A sphere has no overlap of positive surface measure with a
nontrivial translate. This fails for the one-dimensional interval:
a translation can send one endpoint to the other. It also fails for
opposite sides of a rectangle. A boundary crossing can then exchange
spectral components, and the contradiction from strict inclusion does not follow.
\end{remark}

\begin{remark}[A more general boundary condition]\label{rem:general-boundary}
Let $\Omega\subset\R^n$ be a nonempty bounded open set with
$|\partial\Omega|=0$. Suppose its boundary carries a finite nonzero
positive Borel measure $\nu$ such that
\[
 \nu\bigl(\partial\Omega\cap(\partial\Omega-\theta)\bigr)=0
 \qquad(\theta\ne0).
\]
Assume also that, for $\nu$-almost every $t\in\partial\Omega$, every
neighborhood of $t$ meets both $\Omega$ and
$\R^n\setminus\overline\Omega$ in sets of positive Lebesgue measure.
Then $L^2(\Omega)$ admits no Riesz basis of exponentials.
Indeed, Sections~\ref{sec:bumps}--\ref{sec:distance} apply under these
hypotheses. Tonelli's theorem, applied with $\nu$, supplies a point satisfying both
the spectral boundary-null condition and the stated two-sided
neighborhood condition. The latter provides the two sequences of good
parameters. Their symbols
have exactly the limits \eqref{eq:mminus}--\eqref{eq:mplus}, and the same
arguments using strong limits and norm continuity give a contradiction.
In particular one may use surface measure whenever it is finite and
nonzero and satisfies the displayed overlap condition and the two-sided neighborhood condition.
\end{remark}

\section{Triangles}\label{sec:triangles}

\lead{For a triangle, tangential translations can produce overlap of
positive length along an edge, so the sphere argument does not apply
verbatim. We instead isolate one edge from the other two. Crossing it
then brings all the surviving spectral boundary mass into the cutoff
at once, again giving strictly nested limiting ranges.}

\begin{theorem}\label{thm:triangle}
No nondegenerate triangle in $\R^2$ admits a Riesz basis of exponentials.
\end{theorem}

\begin{proof}
An invertible affine change of variables preserves the Riesz basis
property. Indeed, the pullback under $x=a+Ay$ sends
$e_\lambda$ to $e^{2\pi i\lambda\cdot a}e_{A^{\mathsf T}\lambda}$
and multiplies the $L^2$ norm by $|\det A|^{-1/2}$. The phase factors
have modulus one, and the frequency map is invertible. Thus we may assume
\[
 \Omega=\{x_1>0,\ x_2>0,\ x_1+x_2<1\}.
\]
Suppose it has an exponential Riesz basis. Since $\Omega$ is bounded
and its boundary has area zero, Sections~\ref{sec:bumps}--\ref{sec:distance}
give the control measure $\sigma$, the conull set $\mathcal G$,
and the projections with \eqref{eq:stationary-distance-bound}.

\lead{We first choose a point in the interior of the horizontal edge such
that the other two edges carry no translated spectral mass. Their
nonparallel supporting lines make the required Tonelli argument possible.}

Let $L=\{(a,0):0<a<1\}$ and $F=\partial\Omega\setminus L$.
The set $F$ consists of the two other closed edges, including all
vertices. For $\nu=\mathcal H^1|_L$ and every $\theta\in\R^2$,
\[
 \nu\bigl(L\cap(F-\theta)\bigr)=0,
\]
because a horizontal segment meets a translated vertical or diagonal
supporting line in at most one point. Tonelli's theorem therefore gives
\[
 \int_L\sigma(F-t)\,d\nu(t)
 =\int_{\R^2}\nu\bigl(L\cap(F-\theta)\bigr)\,d\sigma(\theta)=0.
\]
Choose $t_0\in L$ with $\sigma(F-t_0)=0$. This set is null for every
$\sigma_f$ by \eqref{eq:control-null}.

\lead{Once the other edges have been removed as spectral null sets, the
sign of one normal coordinate controls every surviving boundary point.
Good parameters approaching from opposite sides therefore have ordered
limiting symbols.}

Near $t_0$, the triangle coincides with the half-plane $x_2>0$.
Choose good parameters $t_j^+\in\mathcal G\cap\Omega$ and
$t_j^-\in\mathcal G\cap(\R^2\setminus\overline\Omega)$ tending to
$t_0$ inside such a neighborhood. Their second coordinates are positive
and negative, respectively. For each point $t_0+\theta\in L$, membership
in the triangle is likewise determined locally by positivity of the second
coordinate. Therefore
$\one_\Omega(t_j^++\theta)\to1$ and
$\one_\Omega(t_j^-+\theta)\to0$ on $L-t_0$.
Away from the boundary both symbols converge to
$\one_\Omega(t_0+\theta)$. All other boundary points belong to the
common null set $F-t_0$.

The limiting symbols are consequently
\[
 m_-=\one_{\Omega-t_0},\qquad
 m_+=\one_{\Omega-t_0}+\one_{L-t_0}.
\]
They are indicators of nested sets.

\lead{The remaining steps are the same as for the ball: the difference
formula gives strong limits, the constant vector makes their ranges
strictly nested, and norm continuity supplies a common nearby projection.}

Formula~\eqref{eq:Pi-difference}
and dominated convergence make $(\Pi_{t_j^\pm}f)_j$ Cauchy for every
$f$. As in Section~\ref{sec:sphere}, the resulting strong limits
$R_-,R_+$ are orthogonal projections with
\[
 \|R_\pm f\|^2=\int m_\pm\,d\sigma_f.
\]
Thus $\|R_-f\|\leq\|R_+f\|$ for every $f$, which gives
$\Ran R_-\subset\Ran R_+$ by taking kernels and orthogonal complements.
Since $0\in L-t_0$ and $\sigma_{\one}=\delta_0$, we have
$R_-\one=0$ and $R_+\one=\one$; the inclusion is strict.
Finally $M_{t_j^\pm}\to M_{t_0}$ in norm, so
\eqref{eq:stationary-distance-bound} gives
$\|R_\pm-M_{t_0}\|\leq\gamma<1$, contradicting
Lemma~\ref{lem:inclusion}.

In this case, $\|(R_+-R_-)f\|^2=\sigma_f(L-t_0)$. The jump may
contain nonzero tangential frequencies, but the constant vector still
guarantees that the jump is nonzero.
\end{proof}

\begin{remark}[Convex polygons with an unpaired edge]\label{rem:odd-polygons}
The same proof excludes an exponential Riesz basis on the interior of a
convex polygon having an edge with no parallel partner. Here edges are
maximal boundary segments. In particular the conclusion holds for
every convex polygon with an odd number of sides: a bounded convex polygon
has at most two maximal edges in each unoriented direction, so an odd
number of sides cannot be partitioned into parallel pairs.

To verify the extension, let $L$ be the relative interior of an unpaired
edge and put $F=\partial\Omega\setminus L$. Every other edge is
nonparallel to $L$, so for $\nu=\mathcal H^1|_L$,
\[
 \nu\bigl(L\cap(F-\theta)\bigr)=0\qquad(\theta\in\R^2).
\]
The intersections with the translated supporting lines are finite;
the two endpoints contribute only finitely many further points.
Tonelli's theorem gives $t_0\in L$ with $\sigma(F-t_0)=0$.
Fix the inward unit normal to $L$. At each point of $L$, membership in
$\Omega$ is locally determined by positive displacement
in this normal direction. Choose good parameters approaching $t_0$
from the two local half-planes. For every $t_0+\theta\in L$, the same
normal sign determines their eventual membership. The symbols are
therefore the same $m_-,m_+$ as in the triangle proof, with this edge $L$.
Their limiting projections satisfy $\Ran R_-\subsetneq\Ran R_+$
because $0\in L-t_0$.
Norm continuity of $M_t$ again supplies the common comparison projection
$M_{t_0}$ and Lemma~\ref{lem:inclusion} gives the contradiction.
\end{remark}

\section{Bounded convex domains with \texorpdfstring{$C^2$}{C2} boundary}\label{sec:convex}

\lead{The boundary criterion extends the nonexistence result from balls to
all nonempty bounded convex open sets with $C^2$ boundary. The key is to
use surface measure on a small positively curved patch, rather than on
the whole boundary, which may contain flat pieces. We first find such a
patch and then show that its intersection with every nontrivial translate
of the boundary has surface measure zero.}

\begin{corollary}\label{cor:convex}
Let $n\geq2$ and let $\Omega\subset\R^n$ be a nonempty bounded convex
open set with $C^2$ boundary. Then $L^2(\Omega)$ admits no Riesz basis of
exponentials.
\end{corollary}

Neither strict convexity of the whole domain nor positive curvature at
every boundary point is required.

\lead{A ball containing the closure of the domain touches its boundary
at a point where all principal curvatures are positive. Continuity gives
a patch with the same property. Convexity then makes each supporting
hyperplane along this patch meet the closed domain at only one point.}

\begin{lemma}\label{lem:convex-patch}
Under the hypotheses of Corollary~\ref{cor:convex}, put
$K=\overline\Omega$ and $S=\partial\Omega$. There is a nonempty relatively
open patch $U\subset S$ on which the second fundamental form is positive
definite, with the convention that convex principal curvatures are
nonnegative. Moreover, if $N(x)$ denotes the outward unit normal, then
\begin{equation}\label{eq:convex-support}
 K\cap\{z:N(x)\cdot(z-x)=0\}=\{x\}
 \qquad(x\in U).
\end{equation}
\end{lemma}

\begin{proof}
Choose $c\in\Omega$ and let
\[
 R=\max_{z\in K}|z-c|>0.
\]
Then $K\subset\overline B(c,R)$, with contact at some $p\in S$.
The tangent hyperplane to the ball at $p$ supports $K$ and, since $S$
is $C^2$, is also the tangent hyperplane to $S$. After a translation
and rotation, we may assume that
\[
 p=0,\qquad c=-Re_n,\qquad N(p)=e_n.
\]
Locally write $S$ as $x_n=f(x')$, with $f\in C^2$,
$f(0)=0$ and $Df(0)=0$. Containment in the ball gives
\[
 f(x')\leq -R+\sqrt{R^2-|x'|^2},
\]
with equality at $x'=0$. Taking second derivatives at this contact point
in the sense of quadratic forms yields
\[
 D^2f(0)\leq-\frac1R I.
\]
Thus the second fundamental form at $p$, represented in these coordinates
by $-D^2f(0)$, is positive definite. By continuity it remains positive
definite on a sufficiently small relatively open patch $U$.

To prove \eqref{eq:convex-support}, fix $x\in U$ and suppose that another
point $z\in K$ lies in the supporting hyperplane at $x$. Convexity places
the segment $[x,z]$ in $K$ and in that hyperplane, hence in $S$.
In coordinates centered at $x$ with horizontal tangent hyperplane, the
boundary is the graph of a function $g$ with $g(0)=0$, $Dg(0)=0$ and
$D^2g(0)$ negative definite. The segment would give
$g(tv)=0$ for a nonzero tangent vector $v$ and all sufficiently small
$t>0$, contradicting
\[
 g(tv)=\frac{t^2}{2}\,D^2g(0)[v,v]+o(t^2)<0.
\]
This proves the singleton supporting-face property.
\end{proof}

\lead{We now compare the outward normals at two boundary points separated
by a fixed translation. Nonparallel normals give a transverse
intersection of codimension two. Equal normals are ruled out by the
singleton supporting-face property, while opposite normals force the
entire intersection to consist of one point.}

\begin{lemma}\label{lem:convex-translates}
Let $U$ be the patch in Lemma~\ref{lem:convex-patch}. Then
\begin{equation}\label{eq:convex-intersection}
 \mathcal H^{n-1}\bigl(U\cap(S-\theta)\bigr)=0
 \qquad(\theta\in\R^n\setminus\{0\}).
\end{equation}
\end{lemma}

\begin{proof}
Fix $\theta\ne0$. At any $x\in U\cap(S-\theta)$, put $y=x+\theta\in S$.
The normals to $S$ and $S-\theta$ at $x$ are $N(x)$ and $N(y)$,
respectively.

If these normals are not parallel, the implicit function theorem makes
the intersection locally a $C^1$ submanifold of dimension $n-2$.
The set of all such transverse intersection points is covered by
countably many of these neighborhoods, since $\R^n$ is second countable.
It therefore has zero $\mathcal H^{n-1}$-measure.

The equality $N(y)=N(x)$ is impossible. Indeed, $x$ and $y$ would
maximize the same linear functional on $K$ and so lie in the same
supporting hyperplane, contrary to \eqref{eq:convex-support} and $x\ne y$.

Finally, suppose $N(y)=-N(x)$ and write $v=N(x)$. The supporting
inequalities at $x$ and $y$ give
\[
 K\subset\{z:v\cdot(z-x)\leq0\},\qquad
 K-\theta\subset\{z:v\cdot(z-x)\geq0\}.
\]
Consequently, \eqref{eq:convex-support} implies
\[
 K\cap(K-\theta)
 \subset K\cap\{z:v\cdot(z-x)=0\}
 =\{x\}.
\]
Since $x$ belongs to both boundaries, $S\cap(S-\theta)=\{x\}$ in this
case. A singleton has zero $\mathcal H^{n-1}$-measure because $n\geq2$.
If this case does not occur, every intersection point in $U$ is
transverse. This proves \eqref{eq:convex-intersection}.
\end{proof}

\lead{Surface measure restricted to the patch now satisfies the
non-overlap hypothesis of Remark~\ref{rem:general-boundary}. The remaining
hypotheses follow from boundary regularity and convexity, so the analytic
construction and the final projection contradiction need no changes.}

\begin{proof}[Proof of Corollary~\ref{cor:convex}]
The compact $C^2$ hypersurface $S$ has zero $n$-dimensional Lebesgue
measure and finite $(n-1)$-dimensional surface measure. Since $U$ is a
nonempty relatively open patch, the positive Borel measure
\[
 \nu=\mathcal H^{n-1}|_U
\]
is finite and nonzero. Lemma~\ref{lem:convex-translates} gives
\[
 \nu\bigl(S\cap(S-\theta)\bigr)
 =\mathcal H^{n-1}\bigl(U\cap(S-\theta)\bigr)=0
 \qquad(\theta\ne0).
\]
Every neighborhood of a point $t\in S$ meets $\Omega$ in a nonempty
open set, because $t$ is a boundary point. It also contains a nonempty
open subset of the strict exterior: the supporting inequality
$N(t)\cdot(z-t)\leq0$ holds throughout $K$, so points sufficiently close
to $t$ on the other side of this hyperplane lie outside $K$.
Both intersections therefore have positive Lebesgue measure.
All the hypotheses of Remark~\ref{rem:general-boundary} hold, and its
conclusion proves the corollary.
\end{proof}

\section{Complete and minimal exponentials on bounded convex domains}\label{sec:complete-minimal}

For domains already known to admit exponential Riesz bases, existence
of a complete minimal system follows at once; this includes the convex
polytopes with centrally symmetric faces treated in
\cite{DebernardiLev}. The theorem below imposes no such geometric
restriction.

\subsection{Statement and Fourier conventions}
A family of function in $L^2(\Omega)$ is \emph{minimal} when
no member belongs to the closed span of the others. Equivalently,
it admits individual biorthogonal functions; no uniform bound on
their norms is required.

\begin{theorem}\label{cm:thm:main}
Let $n\geq1$ and let $\Omega\subset\R^n$ be nonempty, bounded,
open, and convex. There is a locally finite sequence $\Lambda\subset\R^n$ such that
$\{e_\lambda:\lambda\in\Lambda\}$ is complete and minimal in
$L^2(\Omega)$. Its biorthogonal functions are supported in
the John ellipsoid of $\overline\Omega$.
\end{theorem}

The statement concerns ordinary minimality. It does not assert a
Riesz basis, a Schauder basis, or completeness of the biorthogonal
family in $L^2(\Omega)$.

For $n=1$, the domain is an interval $(a,b)$ and the frequencies
$(b-a)^{-1}\mathbb Z$ give an orthogonal basis. We therefore assume
$n\geq2$ throughout the construction below. No smoothness or
strict convexity of the boundary is assumed.

We use the Fourier convention
\[
 \widehat f(z)=\int_{\R^n}f(x)e^{-2\pi ix\cdot z}\,dx,
 \qquad z\in\C^n,
\]
for compactly supported integrable functions, and its distributional
extension. Two classical Fourier facts will be used explicitly.
First, an entire function on $\C^n$ of finite exponential type
which is polynomially bounded on $\R^n$ is the Fourier transform
of a compactly supported distribution. This is the
Paley--Wiener--Schwartz theorem; this formulation is recalled in
\cite[Introduction, pp.~1--2]{Wiegerinck}.
Second, for nonzero compactly supported distributions $u,v$,
\begin{equation}\label{cm:eq:titchmarsh}
 \conv\supp(u*v)=\conv\supp u+\conv\supp v.
\end{equation}
This is the Titchmarsh--Lions convolution support theorem; see
\cite[Theorem 4.2]{Hormander}.
Standard Bessel formulas are taken from \cite{DLMF}.

\Needspace{10\baselineskip}
\subsection{Normalization by the John ellipsoid}

Write $B=\{x\in\R^n:|x|<1\}$ and $v_n=|B|$.
By John's theorem, the convex body $\overline\Omega$ has a unique
inscribed ellipsoid of maximal volume, its \emph{John ellipsoid}; see
\cite[Theorem~3.1, p.~13]{BallConvexGeometry}.
Normalizing this ellipsoid to $\overline B$, uniqueness gives the
following statement.

\begin{lemma}\label{cm:lem:ellipse}
After an invertible affine change of variables, $\overline B$ is the
John ellipsoid of $\overline\Omega$ and
\begin{equation}\label{cm:eq:rigidity}
 B\subset\Omega,
 \qquad
 \{t\in\R^n:\overline B+t\subset\overline\Omega\}=\{0\}.
\end{equation}
\end{lemma}

We henceforth work in this normalization. Denote by $PW_\Omega$
the Fourier image of $L^2(\Omega)$, with physical-space functions
extended by zero. Every member is entire of finite exponential type.
It and its first derivatives are bounded on $\R^n$, since
$\Omega$ is bounded and $L^2(\Omega)\subset L^1(\Omega)$.

\subsection{Spheres and finite-dimensional interpolation spaces}

Set $\nu=n/2-1$ and $\beta=(n+1)/2$. Let $a_j$, $j\geq1$,
be the increasing positive zeros of $J_{\nu+1}=J_{n/2}$.
They are simple, tend to infinity, and $J_\nu(a_j)\ne0$.
The Fourier transform of the ball gives
\begin{equation}\label{cm:eq:definitions}
 \begin{aligned}
 \kappa(z)&=\frac{\widehat{\one_B}(z)}{v_n},
 &\kappa(\xi)&=k(2\pi|\xi|),\\
 k(s)&=2^{\nu+1}\Gamma(\nu+2)\frac{J_{\nu+1}(s)}{s^{\nu+1}},
 &q(z)&=4\pi^2\sum_{k=1}^n z_k^2.
 \end{aligned}
\end{equation}
Here $k(0)=1$, and $\cP_m$ denotes complex polynomials of total
degree at most $m$ in $n$ variables. The Bessel quotient defining
$k$ is understood by its even power series; thus $k$ is entire,
also when $n$ is odd. The same series defines $\kappa$ on $\C^n$
without choosing a square root. In $q$ the squares are complex
squares, with no conjugation. Define
\[
 Q_N(z)=\prod_{j=1}^N(q(z)-a_j^2),\qquad Q_0=1,
 \qquad S_j=\{\xi\in\R^n:2\pi|\xi|=a_j\},
\]
and $V_N=\{P\kappa/Q_N:P\in\cP_{2N}\}$.
The Bessel product gives
\begin{equation}\label{cm:eq:product}
 \kappa(z)=\prod_{j=1}^\infty\left(1-\frac{q(z)}{a_j^2}\right),
\end{equation}
locally uniformly on $\C^n$.

\Needspace{8\baselineskip}
\begin{lemma}\label{cm:lem:membership}
The spaces $V_N$ are contained in $PW_B\subset PW_\Omega$.
They are nested, and
\[
 \dim V_N=\dim\cP_{2N}=\binom{2N+n}{n}.
\]
The function $\kappa/Q_N$ is nonzero at every point of
$\{0\}\cup\bigcup_{j\leq N}S_j$, and every member of $V_N$
vanishes on $S_j$ for $j>N$.
\end{lemma}

\begin{proof}
For a positive zero $a$ of $J_{\nu+1}$ define
\[
 \phi_a(x)=\frac1{v_n a^2}
       \left(\frac{|x|^{-\nu}J_\nu(a|x|)}{J_\nu(a)}-1\right)\one_B(x).
\]
The value at the origin is defined by continuity; the interior
expression is smooth there by its power series in $|x|^2$.
Both its value and its normal derivative vanish on $\partial B$, since
\[
 \frac{d}{dr}\bigl(r^{-\nu}J_\nu(ar)\bigr)
       =-a r^{-\nu}J_{\nu+1}(ar).
\]
The Bessel equation therefore gives, including the zero extension,
\[
 (-\Delta-a^2)\phi_a=\one_B/v_n.
\]
There are no boundary distributions in this equation. Fourier
transformation shows that
$\widehat{\phi_a}=\kappa/(q-a^2)$, with removable singularities.
Partial fractions in $q$ now show that $\kappa/Q_N$ has a
compactly supported inverse transform in $\overline B$.

The real-axis Bessel asymptotic gives
\begin{equation}\label{cm:eq:decay}
 \kappa(\xi)/Q_N(\xi)=O(|\xi|^{-2N-\beta})
                  \quad(|\xi|\longrightarrow\infty).
\end{equation}
Thus $P\kappa/Q_N\in L^2(\R^n)$ for $\deg P\leq2N$:
the radial majorant for its squared norm at infinity is
$r^{n-1-2\beta}=r^{-2}$.
Its inverse distribution is obtained by differentiating a
distribution supported in $\overline B$, hence has the same support
inclusion. Plancherel makes it an $L^2$ function supported there.
Since $\partial B$ has Lebesgue measure zero, this proves membership
in $PW_B$, including $N=0$ directly.

Nesting follows by multiplying the numerator by $q-a_{N+1}^2$.
Multiplication by $\kappa/Q_N$, nonzero on an open set, is injective
on polynomials. The dimension and the remaining assertions follow
from the simple zeros and \eqref{cm:eq:product}.
\end{proof}

The spaces $V_N$ form a nested finite-dimensional interpolation
hierarchy, analogous to a multiresolution decomposition. Here
each level adds polynomial  data on one further Bessel-zero
sphere.

\subsection{Division and the support of the quotient}

We first pass from real spheres to their complex zero hypersurfaces.
This step is valid in every dimension $n\geq2$.

\begin{lemma}\label{cm:lem:quadric}
Let $r>0$. If an entire function $F$ on $\C^n$ vanishes on
$\{\xi\in\R^n:|\xi|=r\}$, then it vanishes on the complex quadric
\[
 \mathcal Q_r=\left\{z\in\C^n:\sum_{k=1}^n z_k^2=r^2\right\}.
\]
It is locally divisible by $\sum_{k=1}^n z_k^2-r^2$ at every point
of this quadric.
\end{lemma}

\begin{proof}
For $z=x+iy\in\mathcal Q_r$, with $x,y\in\R^n$, one has
$x\cdot y=0$ and $|x|^2-|y|^2=r^2$. In particular $x\ne0$.
There is a real two-dimensional linear subspace $E$ containing
$x$ and $y$; if $y=0$, choose any such subspace containing $x$.
In orthonormal coordinates on its complexification, the quadric
section is parametrized by
\[
 z_1=\frac r2(w+w^{-1}),\qquad
 z_2=\frac r{2i}(w-w^{-1}),\qquad w\in\C\setminus\{0\}.
\]
This parametrization is onto: its inverse is $w=(z_1+iz_2)/r$,
which is nonzero because $(z_1+iz_2)(z_1-iz_2)=r^2$.
The unit circle in the $w$ plane maps onto the real circle
$E\cap rS^{n-1}$. The one-variable identity theorem therefore
shows that $F$ vanishes on the whole section, in particular at $z$.
Finally, the gradient of the defining polynomial is $2z\ne0$ on
$\mathcal Q_r$. Local holomorphic coordinates make that polynomial
one coordinate; vanishing when this coordinate is zero gives the
asserted local divisibility.
\end{proof}

We next prove the growth estimates needed for division. Entireness
of a quotient alone does not give compact support of its inverse
transform.

\begin{lemma}\label{cm:lem:type}
Suppose $A,D$ are entire on $\C^n$,
$|A(z)|+|D(z)|\leq C e^{R|z|}$, and $D(0)\ne0$.
If $A/D$ extends to an entire function, that quotient has finite
exponential type.
\end{lemma}

\begin{proof}
Divide both $A$ and $D$ by $D(0)$, adjusting constants, so $D(0)=1$.
For a complex unit vector $u$, restrict the functions to the line
$w\mapsto wu$. Jensen's formula gives
\[
 \frac1{2\pi}\int_0^{2\pi}\log|D(re^{it}u)|\,dt\geq0.
\]
Hence the mean of $\log^-|D|$ on that circle is at most the mean
of $\log^+|D|$. Writing $G=A/D$, we obtain, uniformly in $u$,
\[
 \frac1{2\pi}\int_0^{2\pi}\log^+|G(re^{it}u)|\,dt
       \leq 2\log^+C+2Rr.
\]
The Poisson inequality for the subharmonic function
$\log|G(wu)|$, using a circle of radius $2|w|$ for $|w|\geq1$,
bounds its value by three times the last mean. This gives
$|G(z)|\leq C' e^{R'|z|}$ for $|z|\geq1$, uniformly in the
direction $u$; the compact remainder is harmless. Radii passing
through zeros are handled by limits.
\end{proof}

\begin{lemma}\label{cm:lem:division}
If $F\in PW_\Omega$ vanishes on every $S_j$, $j>N$, then
\begin{equation}\label{cm:eq:quotient}
 G=Q_NF/\kappa
\end{equation}
extends to an entire function which is the Fourier transform of a
compactly supported distribution.
\end{lemma}

\begin{proof}
By Lemma~\ref{cm:lem:quadric}, the numerator $Q_NF$ vanishes on all
the complex quadrics $\mathcal Q_r$ for
$r=a_j/(2\pi)$: for $j\leq N$ the corresponding factor is in
$Q_N$, and for $j>N$ use the hypothesis.
They are exactly the zero hypersurfaces of $\kappa$ in
\eqref{cm:eq:product}, each of multiplicity one, and they do not
intersect each other. They form a locally finite collection because
$a_j\longrightarrow\infty$. Local divisibility therefore removes
all singularities of \eqref{cm:eq:quotient}.

The numerator $Q_NF$ and the denominator $\kappa$ have finite
exponential type. Since $\kappa(0)=1$, Lemma~\ref{cm:lem:type}
gives the same property for $G$. We next prove that $G$ is
polynomially bounded on $\R^n$.

The standard Bessel asymptotics and their differentiated versions
give, for sufficiently large $s$,
\begin{equation}\label{cm:eq:bessellower}
 |k(s)|\geq c(1+s)^{-\beta}
           \min\{1,\dist(s,\{a_j:j\geq1\})\}.
\end{equation}
For completeness, with $c_n>0$ and $\gamma_n=\pi(n+1)/4$,
\[
 \begin{aligned}
 k(s)&=c_n s^{-\beta}\bigl(\cos(s-\gamma_n)+O(s^{-1})\bigr),\\
 k'(s)&=c_n s^{-\beta}\bigl(-\sin(s-\gamma_n)+O(s^{-1})\bigr).
 \end{aligned}
\]
The zeros satisfy $a_j=(j+(n-1)/4)\pi+O(j^{-1})$.
Off fixed small neighborhoods of these zeros the cosine is bounded
away from zero. Inside each neighborhood the derivative has constant
sign and magnitude at least $c s^{-\beta}$; integration from the
zero proves \eqref{cm:eq:bessellower}.

For $\theta\in S^{n-1}$, set
\[
 H_\theta(s)=\prod_{j\leq N}(s^2-a_j^2)\,
                       F(s\theta/(2\pi)).
\]
Uniform boundedness of $F$ and its first real derivatives gives
\[
 |H_\theta(s)|+|H_\theta'(s)|\leq C(1+s)^{2N}
                      \qquad(s\geq0),
\]
uniformly in $\theta$. Also $H_\theta(a_j)=0$ for every $j$.
Near a zero, integrate this derivative estimate from that zero;
away from the zeros use the bound on $H_\theta$ itself. Combining
with \eqref{cm:eq:bessellower} proves
\[
 |G(\xi)|\leq C(1+|\xi|)^{2N+\beta}\qquad(\xi\in\R^n).
\]
The bounded real-frequency region is covered by continuity.
The stated Paley--Wiener--Schwartz theorem now gives a compactly
supported distribution with Fourier transform $G$.
\end{proof}

\begin{proposition}\label{cm:prop:tail}
For a domain satisfying \eqref{cm:eq:rigidity},
\begin{equation}\label{cm:eq:tail}
 \{F\in PW_\Omega:F|_{S_j}=0\text{ for all }j>N\}=V_N.
\end{equation}
\end{proposition}

\begin{proof}
Let $F=\widehat f$ belong to the left side, and suppose $F\ne0$.
By Lemma~\ref{cm:lem:division}, $G=Q_NF/\kappa=\widehat g$ for
a nonzero compactly supported distribution $g$. Put $u=\one_B/v_n$.
The product identity $\kappa G=Q_NF$ gives
\begin{equation}\label{cm:eq:convolution}
 u*g=\prod_{j\leq N}(-\Delta-a_j^2)f.
\end{equation}
The right side is supported in $\overline\Omega$.
The convolution support theorem \eqref{cm:eq:titchmarsh} implies
\[
 \overline B+\conv\supp g
      =\conv\supp(u*g)\subset\overline\Omega.
\]
Condition \eqref{cm:eq:rigidity} forces $\supp g\subset\{0\}$.
A distribution supported at one point is a finite linear combination
of derivatives of its Dirac measure. Hence $G=P$ is a polynomial,
and
\[
 F=P\kappa/Q_N.
\]

We check the degree rather than assuming it from support. If
$m=\deg P$, its leading homogeneous part has absolute value bounded
below on some open patch of $S^{n-1}$. Uniformly on a smaller patch,
$|P(r\theta)|\geq c r^m$ for large $r$.
On real radial intervals of fixed positive length repeated at
bounded spacing, the Bessel asymptotic gives
$|\kappa(r\theta)|\geq c r^{-\beta}$.
Also $|Q_N(r\theta)|\leq C r^{2N}$ for large $r$.
Thus the polar integral of $|F|^2$ on this patch and these radial
intervals is bounded below by a constant times the integral of
\[
 r^{2m-4N-2\beta+n-1}=r^{2m-4N-2}.
\]
If $m\geq2N+1$, this diverges. Since $F\in L^2(\R^n)$, it
follows that $m\leq2N$. In particular, the degree bound is the
same in every dimension.
Therefore $F\in V_N$. The zero function is immediate, and the
reverse inclusion is Lemma~\ref{cm:lem:membership}.
\end{proof}

\subsection{Selection of the frequencies}

Call a finite set
\[
 A\subset\{0\}\cup\bigcup_{j\leq m}S_j
\]
an interpolation set at level $m$ if evaluation on $A$ maps
$\cP_{2m}$ isomorphically onto $\C^A$. By
Lemma~\ref{cm:lem:membership}, this is equivalent to interpolation
by $V_m$. At level zero take $A=\{0\}$.

\begin{lemma}\label{cm:lem:available}
For such a set $A$, the exponentials with frequencies in
$A\cup\bigcup_{j>m}S_j$ have dense span in $L^2(\Omega)$.
\end{lemma}

\begin{proof}
The Fourier transform of a function orthogonal to all those
exponentials belongs to $V_m$ by Proposition~\ref{cm:prop:tail}.
Its vanishing on $A$ and interpolation then force it to be zero.
\end{proof}

\begin{lemma}\label{cm:lem:extension}
Let $A$ be an interpolation set at level $m$, and let $C$ be a finite
set of distinct points on spheres with index greater than $m$.
There are $N>m$ and an interpolation set $A'$ at level $N$ such that
\[
 A\cup C\subset A',\qquad
 A'\setminus A\subset\bigcup_{m<j\leq N}S_j.
\]
\end{lemma}

\begin{proof}
Choose $N>m$ containing all sphere indices occurring in $C$ and
satisfying $2N\geq|A\cup C|-1$. Evaluations at these distinct
points are independent on $\cP_{2N}$: for each $\lambda\in A\cup C$,
the polynomial
\[
 \prod_{\substack{\mu\in A\cup C\\\mu\ne\lambda}}
 \frac{(\xi-\mu)\cdot(\lambda-\mu)}{|\lambda-\mu|^2}
\]
has the required Kronecker values and degree $|A\cup C|-1$.

Evaluations on $A\cup\bigcup_{m<j\leq N}S_j$ span the full dual
of $\cP_{2N}$. Indeed, if $P$ vanishes there, successive division
by the sphere equations gives
\[
 P=\prod_{j=m+1}^N(q-a_j^2)R,\qquad\deg R\leq2m.
\]
To justify each division, write $\xi=(\xi',\xi_n)$ and divide in
$\xi_n$ by $|\xi|^2-r^2$. The remainder is
$R_0(\xi')+\xi_nR_1(\xi')$. Its values at
$\xi_n=\pm\sqrt{r^2-|\xi'|^2}$, for $|\xi'|<r$, force both
coefficient polynomials to vanish on an open subset of $\R^{n-1}$,
hence identically. Distinct sphere factors do not vanish on each other, so division
can be repeated. The displayed product is nonzero on $A$, where
$R$ must therefore vanish. Interpolation at level $m$ gives $R=0$.

Extend the independent evaluations on $A\cup C$ to a basis of this
finite-dimensional dual using points of the indicated spheres.
The resulting points constitute $A'$.
\end{proof}

\begin{proof}[Proof of Theorem~\ref{cm:thm:main}]
Work first in the normalization of Lemma~\ref{cm:lem:ellipse}.
Choose a countable dense set $\{f_1,f_2,\ldots\}$ in $L^2(\Omega)$.
Construct strictly increasing integers $N_s$ and nested finite sets
$A_s$, starting from $N_0=0$, $A_0=\{0\}$, so that $A_s$
interpolates at level $N_s$, all new points have sphere index greater
than $N_{s-1}$, and
\begin{equation}\label{cm:eq:approx}
 \dist\bigl(f_\ell,\spanop\{e_\lambda:\lambda\in A_s\}\bigr)
                   <2^{-s}\qquad(1\leq\ell\leq s).
\end{equation}
Given the previous stage, Lemma~\ref{cm:lem:available} supplies a
finite set $C$ of frequencies on later spheres which, together with
$A_{s-1}$, already achieves these finitely many approximations.
Lemma~\ref{cm:lem:extension} then supplies $A_s$ and $N_s$.
Enlargement of the finite span preserves the approximations.

Set $\Lambda=\bigcup_s A_s$. Condition \eqref{cm:eq:approx} proves
completeness. The set is locally finite: the sphere radii tend to
infinity, and spheres with indices $N_{s-1}<j\leq N_s$ are used
only at stage $s$, when finitely many points are selected.

For each $\lambda\in\Lambda$, use the first stage $s$ with
$\lambda\in A_s$. Interpolation by $V_{N_s}$ gives
$F_\lambda\in V_{N_s}$ with
\[
 F_\lambda(\mu)=\delta_{\lambda\mu}\qquad(\mu\in A_s).
\]
Every later frequency lies on a sphere of index greater than $N_s$,
where $F_\lambda$ vanishes identically. Consequently these Kronecker
values hold on all of $\Lambda$. Lemma~\ref{cm:lem:membership}
identifies $F_\lambda$ with the Fourier transform of an $L^2$
function $g_\lambda$ supported in $\overline B$. Extending by zero
inside $\Omega$ gives
\[
 \langle g_\lambda,e_\mu\rangle_{L^2(\Omega)}
                        =\delta_{\lambda\mu},
\]
proving minimality and the asserted support of the biorthogonals.

Finally, if the original domain is $b+A\Omega$ with $A$ invertible,
transport the frequencies to $A^{-\mathsf T}\Lambda$.
The unitary map $f\mapsto |\det A|^{1/2}f(b+A\,\cdot)$ sends
each transported exponential to a nonzero scalar multiple of the
original one. This preserves completeness and minimality, carries
the supports of the biorthogonals into the John ellipsoid of the original domain,
and preserves local finiteness.
\end{proof}

\begin{remark}\label{cm:rem:scope}
The same selected system is complete and minimal on every measurable
set between the inscribed open ellipsoid and $\Omega$, up to null
sets: completeness passes to restrictions by zero extension, and
the biorthogonals supported in the ellipsoid still work. If the ellipsoid
is a proper subset of the convex domain, these biorthogonals cannot
be complete in $L^2(\Omega)$. Completeness of the dual family is not
part of ordinary minimality.
\end{remark}

\end{document}